\documentclass[12pt]{amsart}

\usepackage{amsmath,amsfonts,amsthm}

\input xy

\xyoption{all}
\newtheorem{theo}[section]{Theorem}
\newtheorem{coro}[section]{Corollary}

\newtheorem*{conj}{Crepant Resolution Conjecture}
\newcommand{\la}{\longrightarrow}
\renewcommand{\o}{\mathcal{O}}
\newcommand{\p}{\mathbb{P}}
\renewcommand{\c}{\mathbb{C}}
\newcommand{\z}{\mathbb{Z}}

\renewcommand{\t}{\mathcal{T}}
\newcommand{\inv}[1]{\left\langle{#1}\right\rangle}
\newcommand{\m}[1]{\overline{M}_{0,{#1}}} 

\newcommand{\X}{\mathcal{X}}

\begin{document}
\title{An example of crepant resolution conjecture in two steps}
\author{Renzo Cavalieri}
\address{Department of Mathematics\\
Colorado State University\\
101 Weber Building\\
Fort Collins, CO 80523-1874}
\address{The Mathematical Sciences Research Institute\\
17 Gauss Way\\
Berkeley, CA 94720-5070}
\email{renzo@math.colostate.edu}
\author{Gueorgui Todorov}
\address{Department of Mathematics\\
Princeton University\\
Fine Hall, Washington Road\\
Princeton NJ 08544-1000 USA}
\address{The Mathematical Sciences Research Institute\\
17 Gauss Way\\
Berkeley, CA 94720-5070}
\email{gtodorov@math.princeton.edu}
\date{}
\maketitle 

\begin{abstract}
We study the relation among the genus $0$ Gromov-Witten theories of the three spaces $\mathcal{X}\leftarrow\mathcal{Z}\leftarrow Y$, where
$\mathcal{X}=[\c^2/\z_3]$, $\mathcal{Z}$ is obtained by a weighted blowup at the stacky point of $\mathcal{X}$, and $Y$ is the crepant resolution of the $A_2$ singularity. We formulate and verify a statement similar to the Crepant Resolution Conjecture of Bryan and Graber.   
\end{abstract}

\section*{Introduction}

The Crepant Resolution Conjecture (CRC), losely stating that a Gorenstein orbifold $\mathcal{X}$ and a crepant resolution $Y$ have ``equivalent'' Gromov-Witten theories, has been at the center of attention in the past few years.

This conjecture was originally formulated in physics by Zaslow and Vafa \cite{zaslow,vafa}, and imported to mathematics by Ruan. Bryan and Graber (\cite{brygra}) proposed an elegant formulation, that however holds under the  technical assumption of $\mathcal{X}$ being hard Lefschetz (see \cite[page 2]{brygra}).
Subsequent work of Coates, Corti, Iritani and Tseng (\cite{ccit1,ccit2}) gave rise to a more general statement (well written in \cite{coruan}) in terms of Givental formalism. Several instances of this conjecture have been verified, and a proof of the conjecture for all toric orbifolds was recently announced by Iritani.

This work stems from a somewhat orthogonal question:
\vspace{-0.5cm}
\begin{center}\textit{does the Crepant Resolution Conjecture factor?}\end{center}
\vspace{-0.5cm} More precisely, given a chain of crepant partial resolutions
$$
\mathcal{X}\leftarrow \mathcal{Z}_1\leftarrow \ldots \leftarrow \mathcal{Z}_r\leftarrow Y
$$
is it possible (and meaningful) to obtain Partial Crepant Resolution statements at each step, and the usual CRC as a composition?

We analyze in detail the first non-trivial example of this situation. For the hard Lefschetz orbifold $\mathcal{X}=[\c^2/\z_3]$, with coarse moduli space the $A_2$ surface singularity, the CRC is proven in the strongest possible form in \cite{bgp}: a linear change of variable between the cohomologies and specialization after analytic continuation of the quantum parameters match the Gromov-Witten potentials.

A weighted blowup at the origin of $\mathcal{X}=[\c^2/\z_3]$ produces the total space of the canonical bundle of a weighted projective line with weights $1$ and $2$, which we denote by $\mathcal{Z}$.

\noindent
\textbf{Results}
\begin{enumerate}
	\item \textit{We express in closed form the genus $0$ (extended) Gromov Witten potential $\hat{\mathcal{F}}_\mathcal{Z}$ of $\mathcal{Z}$.}
	\item \textit{We exhibit linear changes of variables and specializations of the quantum parameters that match, after analytic continuation, $\hat{\mathcal{F}}_{Y}$ with $\hat{\mathcal{F}}_\mathcal{Z}$ and $\hat{\mathcal{F}}_\mathcal{Z}$ with $\hat{\mathcal{F}}_\mathcal{X}$
}
\end{enumerate}

The potential is computed using Atyiah-Bott localization to reduce to hyperelliptic Hodge integral computations from \cite{renzogen}.  By imposing the preservation of the Poincare' pairing we produce the change of variables relating the potentials of $Y$ and $\mathcal{Z}$. The extended potential (introduced in \cite[Section 2.2]{brygra}) is used to invert such transformation. The last side of the triangle is then obtained by composition.

\subsection*{Acknowledgements}
The authors would like to thank Arend Bayer, Aaron Bertram, Jim Bryan, Charles Cadman, Y.-P. Lee, Paul Johnson and Rahul Pandharipande for helpful conversations.

\section{Background}
   
In this section we recall the extended Gromov-Witten potential and the Bryan-Graber statement of the CRC. We follow the general notation introduced in \cite{brygra}, and specifically  the notation in \cite[Section 1.1]{bgp}, for the potentials of $[\c^2/\z_3]$ and its resolution.

For an orbifold $\mathcal{X}$, an element $\widehat{\beta}$ of the  orbifold Neron-Severi group $\widehat{NS}_1(\mathcal{X})$  consists of a curve class $\beta=\sum_1^r \beta_iC_i$ in $\mathcal{X}$ and a collection of integers $\{\hat\beta_{1},\ldots, \hat\beta_s\}$, one for each age one component of the inertia stack of $\mathcal{X}$. The class $\widehat{\beta}$ is called effective if the curve is effective and the integers are non-negative. For an effective class $\widehat{\beta}$ the moduli space $\overline{M}_{g,n}(\mathcal{X},\widehat{\beta})$ is parametrizing genus $g$ twisted stable maps to $\mathcal{X}$ with degree $\beta$, $n$ ordered marked points and $\hat{\beta_i}$ unordered twisted points mapping to the $i$-th component of the inertia stack. The extended potential is defined:
\[
\widehat{F}^{\mathcal{X}} (x_{0},\dotsc ,x_{a},q_{1},\dotsc ,q_{r},u_{1},\dotsc ,u_{s}) =\sum
_{n_0,\ldots , n_a=0}^{\infty }\sum _{\widehat\beta } \left\langle
  \delta_0^{n_0}\cdots\delta_a^{n_a} \right\rangle_{\widehat\beta }^{\mathcal{X}}
\frac{x_0^{n_0}}{n_0!}\cdots\frac{x_a^{n_a}}{n_a!} u_{1}^{\hat\beta_1}\dotsb u_{s}^{\hat{\beta}_s} q_1^{\beta_1}\dotsb q_r^{\beta_r}.
\]

If $\delta_1, \ldots,\delta_s$ are a basis for the age one twisted sectors, then: 
$$
\left\langle
  \delta_0^{n_0}\cdots\delta_a^{n_a} \right\rangle_{\widehat\beta }^{\mathcal{X}}=\frac{1}{\prod_1^s \hat{\beta_i}!}\left\langle
  \delta_0^{n_0}, \delta_1^{n_1+\hat\beta_1}, \cdots, \delta_s^{n_s+\hat\beta_s}, , \cdots,  \delta_a^{n_a} \right\rangle_{\beta }^{\mathcal{X}}
$$
i.e. the extended potential is related to the ordinary potential by; 
$$
\widehat{F}^{\mathcal{X}} = F^{\mathcal{X}}(x_0,(x_{1}+u_{1}),\dotsc,(x_{s}+u_s),x_{s+1},\dotsc,x_a,q_1,\dotsc,q_r)
$$

Next we recall the Bryan-Graber version of the Crepant Resolution Conjecture.

\begin{conj}\label{mainconjecture}
Let $\mathcal{X}$ be a Gorenstein orbifold satisfying the hard Lefschetz condition (see \cite[Definition 1.1]{brygra}) and
admitting a crepant resolution $\pi:Y\rightarrow \mathcal{X}$. Assume $Y$ has $r+s$ curve classes, s of which are contracted by the map to $\mathcal{X}$. Denote by $S$ the  set (in an appropriate positive basis of the $H_2(\mathcal{X})$) of contracted curve classes, and by $R$ the set of curve classes that map birationally onto their image.  

There exists a graded linear
isomorphism
\[
L : H^*(Y) \to
  H_{orb}^*(\mathcal{X})
\]
and roots of unity $c_{1}, \ldots , c_s $ such that the
following conditions hold.
\begin{enumerate}
\item The inverse of $L$ extends the map $\pi ^{*}: H^*(\mathcal{X})\to H^{*}
      (Y)$.
\item Regarding the potential function $\widehat{F^Y}$ as a power series in
      $y_{0},\dotsc ,y_{a},q_{1},\dotsc ,q_{r+s}$, the coefficients
      admit analytic continuations from $(q_{1},\dotsc ,q_{s}) =
      (0,\dotsc ,0)$ to $(q_{1},\dotsc ,q_{s}) = (c_{1},\dotsc
      ,c_{s})$.
\item The potential functions $\widehat{F}^\mathcal{X}$ and $\widehat{F}^{Y}$ are equal after
      the substitution
\begin{align*}
y_i &= \sum _{j} L^j_iz_j\\
q_i &= \begin{cases} c_i e^{\sum _{j} L^j_i u_j}  & \text{when } i\in S\\
q_i &\text{when } i\in R. \end{cases}
\end{align*}
\end{enumerate}
\end{conj}

We propose a slightly extended version of this statement that we wish to call the Partial Crepant Resolution Conjecture (\textbf{PCRC}). This differs from the above statement in only two ways:
\begin{enumerate}
	\item We don't require $\pi$ to be a resolution. We wish to apply the statement to any crepant map of Gorenstein orbifolds $\pi:\mathcal{Z}\rightarrow \X$.
	\item We conjecture the need for roots of unity $\tilde{c}_1,\ldots, \tilde{c}_r$ such that the non-special quantum parameters $q_i$ are specialized to $\tilde{c}_i q_i$ in order to match the potentials. 
\end{enumerate}

Let  $\X=\mathcal[\c^2/\z_3]$ where $\z_3$ acts by opposite weights and $Y$ the crepant resolution of the coarse moduli space 
$\c^2/\z_3$. In \cite{bgp}, Bryan, Graber and Pandharipande give a closed form expression for the equivariant genus 0 Gromov-Witten potential of $\mathcal{X}$ (Theorem 1.2),  of $Y$ (Theorem 1.1), and verify the CRC for this pair (Theorem 1.3).

In particular, they exhibit an explicit change of variables that identifies the potentials. Here we tweak this change of variable to the extended potential case: 
\begin{align*}
y_{0}&=x_{0}\\
y_{1}&=\frac{i}{\sqrt{3}} (\omega x_{1}+\overline{\omega} x_{2})\\ 
y_{2}&=\frac{i}{\sqrt{3}} (\overline{\omega}  x_{1}+\omega x_{2})\\
q_{1}&=\omega e^{\frac{i}{\sqrt{3}} (\omega u_{1}+\overline{\omega} u_{2})}\\
q_{2}&=\omega e^{\frac{i}{\sqrt{3}} (\overline{\omega} u_{1}+\omega u_{2})} 
\end{align*}

\section{The Potential}

\subsection{Set-up and Notation}

 Our main computational technique is virtual localization, for which we refer to \cite{grapan}, and moduli spaces of admissible covers, for which we refer to \cite{renzohod}.
Let $\mathcal{Z}$ be the  weighted blow-up of the origin in $\c^2/\z_3$ with weights 
one and two. We consider $\mathcal{Z}$ with the induced orbifold structure. A detailed explanation how this is done in general is in \cite{kerr}. The stack $\mathcal{Z}$ is isomorphic to a local  
$\p(1,2)$: the total space of  the canonical line bundle on the weighted projective line. Let us fix some notations. Let $x_0$ and $x_1$ be homogeneous coordinates on $\p(1,2)$  
so that 
$$
(x_0:x_1)=\left(\frac{x_0}{x_1^{1/2}}:1\right)=\left(1:\frac{x_1}{x_0^2}\right)
$$

 Let $X_0=\frac{x_0}{x_1^{1/2}}$ denote the pseudo-coordinate on the chart around $0$, isomorphic to $[\c/\mathbb{Z}_2]$. Let $X_1=\frac{x_1}{x_0^2}$ 
be the coordinate on the chart around $\infty$, isomorphic to $\c$.
We have the relation 
$
X_0^2=X_1^{-1}.
$

The only stacky point of $\mathcal{Z}$ is at $X_0=0$ and it has   stabilizer group $\z_2$.
When we consider moduli spaces of twisted stable  maps to $\mathcal{Z}$, the only possible non-trivial stabilizer at points on the source curves   is  $\z_2$. Therefore from now on, by ``stacky point'' we 
mean a point with a $\z_2$ stabilizer.

We  call $\o_{\p(1,2)}(1/2)$ the positive generator of 
Pic$\p(1,2)\cong \z$. We denote by $1$ and $H$  the fundamental class and the class of a point in the untwisted sector of $\mathcal{Z}$, by $S$ the fundamental class of the twisted sector. We note that $H$ and $S$ span $H^2_{CR}(\mathcal{Z})$. The dual variables to $1,H,S$ are  denoted $z_0,z_1,z_2$. 

By $\m{n_1,n_2}(X,\beta)$  we denote the components of $\m{n_1+n_2}(X,\beta)$ on which $n_2$ of the evaluation maps go the twisted sector, that is twisted stable maps with $n_1+n_2$ marked points, $n_2$ 
of which are stacky.

The torus action $T=\c^*\times \c^*$  on $\c^2$ with weights $(t_1,t_2)$ induces a natural action on $\mathcal{Z}$. The fixed points for the $T$ action on $\mathcal{Z}$ are $0$ and $\infty$. The weights on the base and fiber directions are: 
$$
\left(t_2-\frac{t_1}{2}, \frac{3t_1}{2}\right)\textrm{ and }\left(t_1-2t_2, 3t_2\right).
$$

% Let us denote by $1$ the fundamental class, and by $S$ and $H$ the stacky and 
%divisorial classes respectively and let $z_0$, $z_1$ and $z_2$ the corresponding 
%dual variables. 
\subsection{Degree $0$ invariants}
\label{zero}
Except for the case with only stacky insertions the degree zero invariants are  given by 
the triple intersections in equivariant cohomology 
$$
\inv{a,b,c}_0=\int_Xa\cup b\cup c,$$
which are computed by localization: 

$$
\inv{1,1,1}_0=\frac{1}{3t_1t_2}
$$
$$
\inv{1,1,H}_0=0=\frac{1}{2}\frac{-t_1}{\tfrac{3t_1}{2}(t_2-\tfrac{t_1}{2})}+\frac{-2t_2}{3t_2(2t_1-2t_2)}
$$
$$
\inv{1,H,H}_0=\frac{1}{2}\frac{(-t_1)^2}{\tfrac{3t_1}{2}(t_2-\tfrac{t_1}{2})}+\frac{(-2t_2)^2}{3t_2(2t_1-2t_2)}=-\frac{2}{3}
$$
$$
\inv{H,H,H}_0=-\frac{2(t_1+2t_2)}{3},\quad \inv{1,S,S}_0=\frac{1}{2},\quad \inv{H,S,S}_0=\frac{-t_1}{2}.
$$

 The most interesting degree zero invariants are given by  $n$-stacky 
insertions. 
Since the insertions are stacky  the corresponding components of the moduli space parametrize maps that contract the curve to $0$. By monodromy considerations (or orbifold Riemann-Roch) the number of insertions $n$ must be even (and in particular $n\ge 4$).
  The invariants are computed  by an integral over the component of 
$\m{n}(B\z_2)$  where all the evaluation maps go to the twisted sector, which we identify with the moduli space of admissible hyperelliptic covers of a rational curve %$Adm_{g\stackrel{2}{\rightarrow} 0,(t_1,\dots,t_{2g+2})}$
(here  $n=2g+2$):  

$$
\inv{S^{2g+2}}_0 = \int_{Adm_{g\stackrel{2}{\rightarrow} 0,(t_1,\dots,t_{2g+2})}}\hspace{0.3cm} e^{eq}(R^1\pi_* f^* (L\oplus L))
$$ 

 The obstruction bundle on admissible covers is obtained by pull-pushing the normal bundle to $0$ (isomorphic to 
%and  $L \oplus L \to B\z_2$ 
two copies of the non-trivial bundle $L$ on $B\z_2$, linearized with weights $\frac{3t_1}{2}$ and $t_2-\frac{t_1}{2}$) via   
% $e(R^1\pi_* f^* (L\oplus L))$
%where $f$ and $\pi$ are 
the universal map and universal curve 
for $\m{n}(B\z_2)$. 
The bundle $R^1\pi_* f^* L$ is isomorphic to the dual of the Hodge bundle pulled back from  the natural forgetful morphism $Adm_{g\stackrel{2}{\rightarrow} 0,(t_1,\dots,t_{2g+2})}\rightarrow \overline{\mathcal{M}}_{g}$. We conclude that 
$$
\inv{S^n}_0=-(t_1+t_2)\int_{Adm_{g\stackrel{d}{\rightarrow} 0,(t_1,\dots,t_{2g+2})}} \lambda_g\lambda_{g-1}.
$$

These integrals  have been computed in \cite[Corollary 2]{fabpan} or \cite{bct}. If we denote this contribution to the potential by $G$, then the first non-zero coefficient is in degree $4$ and:
$$
G'''=\frac{1}{2}\tan\left(\frac{z_2}{2}\right).
$$

\subsection{Positive degree}

 Let $E\cong\p(1,2)$  be the zero section in $\mathcal{Z}$. The image of any non-constant map must lie 
in $E$ and so we have that $\m{n_1,n_2}(X,d[E])\cong \m{n_1,n_2}(\p(1,2),d) $. 
With this identification we have that:

$[\m{n_1,n_2}(X,d[E])]^\textrm{vir}= e(R^\bullet\pi_*f^*N_{E/X})$. The Euler sequence on $\p(1,2)$
%%%$$
%%%0\la \o\p(1,2)\left(-\frac{3}{2}\right)\la \o\p(1,2)\left(-1\right) \oplus \o\p(1,2)\left(-\frac{1}{2}\right)\la \o\p(1,2)\la 0
%%%$$
$$
0\la \o_{\p(1,2)}(-3/2)\la \o_{\p(1,2)}(-1) \oplus \o_{\p(1,2)}(-1/2)\la \o_{\p(1,2)}\la 0
$$
%%%$$
%%%0\la \o(-\frac{3}{2})\la \o(-1) \oplus \o(-\frac{1}{2})\la \o\la 0
%%%$$
 gives the relation
$$
e(R^1\pi_*f^*N_{E/X})=(t_1+t_2)e(R^1\pi_*f^*(\o_{\p(1,2)}(-1) \oplus \o_{\p(1,2)}(-1/2))).
$$

Hence we can express invariants as: 
  $$\inv{H^{n_1}S^{n_2}}_d=(t_1+t_2)\int_{\m{n_1,n_2}(\p(1,2),d)}e(R^\bullet\pi_*f^*(\o_{\p(1,2)}(-1) \oplus \o_{\p(1,2)}(-1/2))).$$
 The advantage of this reduction is that now we are in the local Calabi-Yau case. Because of the divisor equation it is enough to compute the invariants $\inv{S^{n}}_d$, given by the integral:
$$
I_{d,n}:=\int_{\m{0,n}(\p(1,2),d)}e(R^\bullet\pi_*f^*(\o_{\p(1,2)}(-1) \oplus \o_{\p(1,2)}(-1/2))).
$$
%Not set $n_2=n$ and call this integral $I_{d,n}$
% The computation is going to be done by localization.
Consider  an auxiliary $\mathbb{C}^\ast$ action on $\p(1,2)$, inherited from the standard action on the coarse moduli space. In homogeneous coordinates:
 $$
s\cdot(x_0,x_1)\la (s^{1/2}x_0,x_1).
$$
The tangent bundle is canonically linearized with weights 
$1/2$ over $0$ and $-1$ over $\infty$. Table \ref{linearization} shows the choice of the lifting of the action to the bundles $\o_{\p(1,2)}(-1)$ and $\o_{\p(1,2)}(-1/2)$ that we make in order to eliminate many fixed loci. 

\begin{table}[tb]	
	\begin{center}
\begin{tabular}{|l||c|c|}
\hline
      & weight over $0$ & weight over $\infty$ \\
\hline
\hline
 $\o_{\p(1,2)}(-1)$ & 0 & 1 \\
\hline
$\o_{\p(1,2)}(-1/2)$  &  -1/2 & 0  \\
\hline
$T_{\p(1,2)}$  &  +1/2 & 0  \\
\hline
\end{tabular}
\end{center}\caption{Weights of the lifting of the torus action.}
\label{linearization}
\end{table}
 A fixed map for the torus action  consists of a collection of genus zero (orbi-)curves mapping of positive degree to $\p(1,2)$ and fully ramified over $0$ and $\infty$, attached to other genus zero curves contracting to the two fixed points. The marked points must lie on curves contracted over $0$.
Our choice  of weights  forces the possibly contributing localization graphs to have valence one over both $0$ and $\infty$.   The only fixed locus that survives parametrizes maps $f:\mathcal{C} \la \p(1,2)$, where $\mathcal{C}$ has two components $F$ and $C$. The curve $F$ carries all the marked points and it is contracted over zero while $f_{|C}$ is a degree $d$ map fully ramified over $0$ and $\infty$.

%consists of one curve mapping to  $\p(1,2)$ with degree $d$ and a contracted curve over $)$ carrying all the marks. The %description of such a map is slightly different for even and odd degrees, and hence we make two separate discussions.
% have only one edge. {\bf do you have a template for a pic?}

\subsection{Odd degree}
\label{odddeg}

When the degree is odd the node is stacky. Since we must have an even number of stacky points on $\mathcal{C}$, we can only have invariants with an odd number of insertions. The non-contracting curve $C$ is again a weighted projective line with weights one and two. If  we choose  local coordinates $Y_0$ and $Y_1$ the map is given by $f^*X_0=Y_0^d$. The weight of $Y_0$ is $-\frac{s}{2d}$. We describe the localization contribution following \cite{mirr}:
  
 \begin{description}
 	\item[edge terms:] here we have the weights of the $H^1$'s at the numerator, and of the $H^0$'s at the denominator. A basis of $H^1(C,f^*\o(-1/2))=H^1(C,\o(-\frac{d}{2}))$  is given by 
$Y_0^{-2(i+1)}$ for $0\le i \le \frac{d-3}{2}$. With our choice of lifting of weights we obtain:   
$$
\frac{s^{\frac{d-1}{2}}}{(2d)^\frac{d-1}{2}}(d-1)!!.$$
Similarly $H^1(C,f^*\o(-1))=H^1(C,\o(-d))$ is generated by $Y_0^{-2i}$ for $1\le i \le d-1 $. The contribution is:   
$$
\frac{s^{d-1}}{(2d)^{d-1}}(2d-2)!
$$ 
At the denominator $H^0(C,f^*\t_\p)$ is generated by
$Y_0^{2i+1},\textrm{ for }0\le i\le \frac{3d-1}{2}.$
The product of weights is:
$$
-\frac{s^\frac{3d-1}{2}}{(2d)^\frac{3d-1}{2}}(2d)!!(d-1)!!.
$$ 
 	\item[vertex terms:] pulling back the bundles to the contracted component we obtain a trivial bundle on the contracted curve. But the action of $\z_2$ on the fiber of the bundle $f^*\o(-1)$ is trivial, whereas on the other two bundles it is not. After pull-pushing, we obtain two copies of the dual of the Hodge bundle on the genus $g$ hyperelliptic cover corresponding to the map $f_{|E}$, linearized with the appropriate weights. The contribution is therefore the equivariant euler class of
$$
e^{eq}(\mathbb{E}^{\vee}(\tfrac{1}{2})\oplus \mathbb{E}^{\vee}(-\tfrac{1}{2}))=\left(\frac{s}{2}\right)^{2g},
$$
where the equality follows from Mumford's relation (\cite{mumford}).
 	\item[flag terms:] in this case there is no flag contribution since we are pulling back the non-trivial line bundle $L$ on $B\z_2$ via a map with a twisted point.
 	\item[automorphisms and node smoothing:] for the automorphism of the curve we have a term corresponding to moving the vertex over $\infty$. Together with the node smoothing contribution we obtain:
 	$$
 	-\frac{s}{d}\left(\frac{s}{2d}-\frac{\psi}{2}\right)^{-1}
 	$$
 \end{description}

 Putting everything together 
we obtain:
\begin{equation}
I_{d,2g+1}=(-1)^{g+\frac{d-1}{2}}\left(\frac{d}{2}\right)^{2g-1}
\int_{Adm_{g\stackrel{d}{\rightarrow} 0,(t_1,\dots,t_{2g+2})}}\psi^{2g-1}.
\end{equation}
The integral is just the hyperelliptic Hurwitz number with value $\frac{1}{2}$ for any $g$, and hence we have:
\begin{equation}\label{odd}
\mathcal{I}_d=\sum_{g\geq 0}{I}_{d,2g+1}\frac{z_2^{2g+1}}{(2g+1)!}=(-1)^\frac{d-1}{2}\frac{2}{d^3}\sin\left(\frac{dz_2}{2}\right).
\end{equation}
\noindent\textbf{Remark:} the  $g=0$ term in formula (\ref{odd}) corresponds to the invariant with only one stacky insertion, and therefore no contracted curve $F$.

\subsection{Even degree}
\label{even}
In this case $F$ and $C$ are connected by a non-stacky node, and $C$  is a smooth rational curve. If  we choose  local coordinates $Z_0$ and $Z_1$ the map on $C$ is given by $f^*X_0=Z_0^{\frac{d}{2}}$. 
 \begin{description}
 	\item[edge terms:]  we have the weights of the $H^1$'s at the numerator, and of the $H^0$'s at the denominator. %A basis of $H^1(C,f^*\o(-1/2))=H^1(C,\o(-\frac{d}{2}))$  is given by $Y_0^{-2(i+1)}$ for $0\le i \le \frac{d-3}{2}$. 
With our choice of lifting of weights we obtain:   
$$
\frac{s^{\frac{d-1}{2}}}{(2d)^\frac{d-1}{2}}(d-1)!!.
$$
Similarly for $H^1(C,f^*\o(-1))=H^1(C,\o(-d))$ the contribution is:   

$$\frac{s^{d-1}}{d^{d-1}}(d-1)!$$
At the denominator for $H^0(C,f^*\t_\p)$ the product of weights is:
$$
2\frac{s^\frac{3d}{2}}{(2d)^\frac{d}{2}d^d}d!d!!
$$ 
 	\item[vertex terms:] pulling back the bundles to the contracted component we obtain a trivial bundle on the contracted curve. But the action of $\z_2$ on the fiber of the bundle $f^*\o(-1)$ is trivial, whereas on the other two bundles it is not. After pull-pushing, we obtain two copies of the dual of the Hodge bundle on the genus $g$ hyperelliptic cover corresponding to the map $f_{|E}$, linearized with the appropriate weights. The contribution is again:
$$
e^{eq}(\mathbb{E}^{\vee}(\tfrac{1}{2})\oplus \mathbb{E}^{\vee}(-\tfrac{1}{2}))=\left(\frac{s}{2}\right)^{2g}
$$
 	\item[flag terms:] since we are pulling back the non-trivial line bundle $L$ on $B\z_2$ via a map with a regular point we have a flag term:
$$
\frac{s}{2}
$$

 	\item[automorphisms and node smoothing:] for the automorphism of the curve we have a term corresponding to moving the vertex over $\infty$. Together with the node smoothing contribution we obtain:
 	$$
 	-\frac{s}{d}\left(\frac{s}{d}-\psi\right)^{-1}
 	$$
 \end{description}

Since the node is not stacky we also have an additional gluing factor of $2$. Including the invariant with no insertions (corresponding to $g=-1$), we obtain:

$$ 
\mathcal{I}_{d}= \sum_{g=-1}^\infty I_{d,2g+2} \frac{z_2^{2g+2}}{(2g+s)!} =(-1)^{\frac{d}{2}}\frac{2}{d^3}\cos\left(\frac{dz_2}{2}\right).
$$
\subsection{The potential}
Collecting the information from sections \ref{zero}, \ref{odddeg}, \ref{even} and using the divisor equation 
%and the information for the empty invariant which is 
%$$
%\inv{}_d=\frac{1}{2d^3}
%$$
we obtain the following Theorem.

\begin{theo} The equivariant genus 0 Gromov-Witten potential of $\mathcal{Z}$ is 
\begin{eqnarray*}
F^{\mathcal{Z}}&=&\frac{z_0^3}{18t_1t_2}-
\frac{1}{3}z_0z_1^2\\
&&+\frac{1}{4}z_0z_2^2-\frac{t_1}{4}z_1z_2^2-\frac{1(t_1+2t_2)}{9}z_1^3-(t_1+t_2)G\\
&&+(t_1+t_2)\left(\sum_{d\textrm{ odd}}(-1)^\frac{d-1}{2}\frac{2}{d^3}\sin\left(\frac{dz_2}{2}\right)e^{dz_1}q^d+\sum_{d\textrm{ even}}(-1)^{\frac{d}{2}}\frac{2}{d^3}\left(\cos\left(\frac{dz_2}{2}\right)\right)e^{dz_1}q^d \right).
\end{eqnarray*}
where $G'''=\frac{1}{2}\tan(\frac{z_2}{2})$.        
\end{theo}

\section{Checking the PCRC}

\begin{theo}
\label{check}
The change of variables:
\begin{align}
y_0&= z_0 \nonumber\\
y_1&= iz_2 \nonumber\\
y_2&= z_1-\frac{i}{2}z_2\nonumber\\
q_1&=-e^{iu}\nonumber\\
q_2&=iq, \label{cov}
\end{align}
verifies the PCRC for the pair $\mathcal{Z},Y$.
\end{theo}
% if we set $u=0$ we compare the regular (not extended) potentials. 
% We potential of $Y$ after then becomes 
%$$
%\sum_{d=1}^{\infty}\frac{1}{d^3}\left[(-1)^de^{-id(z_2+u)}+e^{-\frac{id}{2}z_2+dz_1}(iq)^d+e^{\frac{id}{2}z_2+iu+dz_1}(-iq)$$

\noindent \emph{Proof.} A straightforward algebraic substitution  checks that (\ref{cov}) matches the degree $0$ part of the potential of $Y$ with the degree $0$ three pointed invariants contribution to the potential of $\mathcal{Z}$, up to a residual contribution of 
\begin{equation}\label{residual}
-\frac{i}{12}(t_1+t_2)z_2^3.
\end{equation}

Next we focus on the term:
$$
(t_{1}+t_{2})\sum _{d=1}^{\infty }\frac{1}{d^{3}}
(e^{y_{1}}q_{2})^{d}
$$

This term appears in the Gromov-Witten potential for $[\mathbb{C}^2/\z_2]$; the change of variables (\ref{cov}) on this term coincides with the change of variables for the case of $[\mathbb{C}^2/\z_2]$. Bryan and Graber prove (\cite[Corollary 3.4]{brygra}) that after the change of variables (and analytic continuation) one obtains:
$$
-(t_1+t_2)G+ \frac{i}{12}(t_1+t_2)z_2^3.
$$

This clears up the term (\ref{residual}) and completes the matching for the degree $0$ invariants of $\mathcal{Z}$.

Plugging (\ref{cov}) into the remaining two summations we obtain: 
$$
\sum_{d=1}^{\infty}\frac{(qe^{z_1})^{d}}{d^3}e^\frac{idu}{2}(i)^d
\left[\frac{e^{-i\frac{dz_2+du}{2}}+(-1)^de^{i\frac{dz_2+du}{2}}}{2} \right]
$$
Hence for $d$ odd the expression in the brackets above becomes 
$$
(-1)^\frac{d-1}{2}\left[\frac{e^{i\frac{dz_2+du}{2}}-e^{-i\frac{dz_2+du}{2}}}{2i} \right]=(-1)^\frac{d-1}{2}\sin \left(\frac{d(z_2+u)}{2}\right)
$$
and for $d$ even: 
$$
(-1)^\frac{d}{2}\left[\frac{e^{i\frac{dz_2+du}{2}}+e^{-i\frac{dz_2+du}{2}}}{2} \right]=(-1)^\frac{d}{2}\cos \left(\frac{d(z_2+u)}{2}\right)
$$
This concludes the proof of Theorem \ref{check}.

From \cite{bgp}, the change of variables between the extended Potentials of $\X$ and $Y$ is given by 

\begin{align}
\label{covbgp}
y_{0}&=x_{0} \nonumber\\
y_{1}&=\frac{i}{\sqrt{3}} (\omega x_{1}+\overline{\omega} x_{2})\nonumber\\ 
y_{2}&=\frac{i}{\sqrt{3}} (\overline{\omega}  x_{1}+\omega x_{2})\nonumber\\
q_{1}&=\omega e^{\frac{i}{\sqrt{3}} (\omega s_{1}+\overline{\omega} s_{2})}\nonumber\\
q_{2}&=\omega e^{\frac{i}{\sqrt{3}} (\overline{\omega} s_{1}+\omega s_{2})} 
\end{align}

One can then obtain the change of variable between $\widehat{F}^{\mathcal{Z}}$ to $\widehat{F}^X$ by inverting (\ref{cov}) and composing with (\ref{covbgp}).
\begin{coro}
The change of variables:
\begin{align*}
z_{0}&=x_{0}\\
z_{1}&=\frac{i}{2\sqrt{3}}\left((\bar\omega-1) x_{1}+(\omega-1) x_{2}\right)\\ 
z_{2}&=\frac{1}{\sqrt{3}}\left({\omega}  x_{1}+\bar\omega x_{2}\right)\\
q&=-i\omega e^{\frac{i}{\sqrt{3}} (\bar\omega s_{1}+{\omega} s_{2})}\\
u&=-\frac{\pi}{3}+\frac{1}{\sqrt{3}} ({\omega} s_{1}+\bar\omega s_{2})
\end{align*}
verifies the PCRC for the pair $\X$ and $\mathcal{Z}$
\end{coro}

\noindent
\textbf{Remark:} These change of variables is not unique, as it depends on a choice for a branch of the logarithm when inverting (\ref{cov}). However no choice  will give $u=0$ when $s_1,s_2=0$, which suggests that the extended potential is not only a natural, but a necessary choice in order to obtain a PCRC statement.

\bibliographystyle{alpha}
\bibliography{ExampleCrepResol11may09}

\end{document}